\documentclass[12pt]{article}

\usepackage{amssymb}

\newcommand{\qed}{$\;\;\;\Box$}
\newenvironment{proof}{\par\smallbreak{\sl Proof.~}}
{\unskip\nobreak\hfill \qed \par\medbreak}

\newcounter{claim}
\renewcommand{\theclaim}{\arabic{claim}}
{\par\medskip\par}

{\qed\par\smallbreak}
\newcommand{\hide}[1]{}

\renewcommand{\L}{{\mbox{L}}}

\newcommand{\N}{{\Bbb N}}
\newcommand{\R}{{\Bbb R}}

\newcommand{\C}{{\Bbb C}}
\newcommand{\Z}{{\Bbb Z}}

\newcommand{\LL}{{\cal L}}
\newcommand{\beq}{\begin{equation}}
\newcommand{\ee}{\end{equation}}

\renewcommand{\d}{\partial}

\newtheorem{thm}{Theorem}
\newtheorem{lemma}[thm]{Lemma}

\newtheorem{rem}[thm]{Remark}

\newcommand{\al}{\alpha}

\newcommand{\ga}{\gamma}
\newcommand{\de}{\delta}
\newcommand{\eps}{\varepsilon}
\newcommand{\vphi}{\varphi}

\newcommand{\om}{\omega}

\newcommand{\im}{\mathop{\rm im}}

\newcommand{\codim}{\mathop{\rm codim}}
\title{
Fredholm Alternative for Periodic-Dirichlet Problems for  Linear
Hyperbolic Systems} 

\newcounter{thesame}
\author{
Irina Kmit \ \ \ Lutz Recke\\
{\small
Institute for Applied Problems of Mechanics and Mathematics,}\\
{\small Ukrainian Academy of Sciences}\\
{\small Naukova St.\ 3b,
79060 Lviv,
Ukraine}
\\
{\small   E-mail:
{\tt kmit@informatik.hu-berlin.de}}\\[5mm]
{\small
Institute of Mathematics, Humboldt University of Berlin,}\\
{\small Rudower Chaussee 25,
12489 Berlin, Germany}\\
{\small   E-mail:
{\tt recke@mathematik.hu-berlin.de}}\\
}

\date{}

\begin{document}

\maketitle

\begin{abstract}
\noindent 
This paper concerns hyperbolic systems of two linear
first-order PDEs in one space dimension with periodicity
conditions in time and reflection boundary conditions in space. 
The coefficients of the PDEs are supposed to be time independent, but 
allowed to be discontinuous with respect to the space variable.
We construct two scales of Banach spaces
(for the solutions and for the right hand sides of the equations, respectively) 
such that
the problem can be modeled by means of Fredholm operators of index zero 
between corresponding spaces of the
two scales.
\end{abstract}

\section{Introduction}\label{sec:intr}
\renewcommand{\theequation}{{\thesection}.\arabic{equation}}
\setcounter{equation}{0}

\subsection{Problem and main results}\label{sec:results}

This paper concerns linear inhomogeneous hyperbolic systems of first order PDEs
in one space dimension of the type
\beq\label{eq:1.1}
\left.
\begin{array}{l}
\partial_tu  + \partial_xu + a(x)u + b(x)v  = f(x,t)\\
\partial_tv  -  \partial_xv + c(x)u + d(x)v  = g(x,t)
\end{array}
\right\}
\; 0 \le x \le 1, \; t \in \R
\ee
with time-periodicity conditions
\beq\label{per}
\left.
\begin{array}{l}
u\left(x,t+\frac{2\pi}{\om}\right) = u(x,t)\\[1mm]
v\left(x,t+\frac{2\pi}{\om}\right) = v(x,t)   
\end{array}
\right\}
\; 0 \le x \le 1, \; t \in \R
\ee
and reflection boundary conditions
\beq\label{bc}
\left.
\begin{array}{l}
u(0,t) = r_0v(0,t)\\
v(1,t) = r_1u(1,t)
\end{array}
\right\}
\; t \in \R.
\ee
Together with the periodic-Dirichlet problem (\ref{eq:1.1})--(\ref{bc}) we consider
its homogeneous adjoint variant
\begin{equation}\label{eq:1.2}
\left.
\begin{array}{l}
-\partial_tu  - \partial_xu + a(x)u + c(x)v  = 0\\
-\partial_tv  + \partial_xv + b(x)u + d(x)v  = 0
\end{array}
\right\}
\;0 \le x \le 1, \; t \in \R,
\ee
\beq \label{adper}
\left.
\begin{array}{l}
u\left(x,t+\frac{2\pi}{\om}\right) = u(x,t)\\[1mm]
v\left(x,t+\frac{2\pi}{\om}\right) = v(x,t)  
\end{array}
\right\}
\;0 \le x \le 1, \; t \in \R,
\ee
\beq \label{adbc}
\left.
\begin{array}{l}
v(0,t) = r_0u(0,t)\\
u(1,t) = r_1v(1,t)
\end{array}
\right\}
\; t \in \R.
\end{equation}
Here
$\om>0$ and $r_0, r_1\in\R$ are fixed numbers, $a,b,c,d:[0,1] \to \R$
are fixed coefficient  functions, and the right-hand sides $f,g:[0,1] 
\times \R \to \R$
are supposed to be $\frac{2\pi}{\om}$-periodic with respect to $t$.

Roughly speaking, we will prove the following: Suppose
\begin{equation}\label{eq:1.4}
|r_0r_1|\ne e^{\int\limits_0^1\left(a(x)+d(x)\right)\,dx}.
\end{equation}
Then there exists a solution to~(\ref{eq:1.1})--(\ref{bc}) if 
and only if the pair
$(f,g)$ of the right hand sides is orthogonal in 
$L^2\left((0,1)\times\left(0,\frac{2\pi}{\om}\right);\R^2\right)$ to 
any solution to~(\ref{eq:1.2})--(\ref{adbc}).
Moreover, the dimension of the space of all solutions to~(\ref{eq:1.1})--(\ref{bc}) 
with
$f=g=0$ is finite and is equal to the dimension of the space of all solutions 
to~(\ref{eq:1.2})--(\ref{adbc}).
More exactly, we construct two scales $V^\gamma(r_0,r_1)$ and $W^\gamma$ (with scale
parameter $\gamma \ge 1$) of Banach spaces such that $V^\gamma(r_0,r_1) \hookrightarrow W^\gamma
 \hookrightarrow C(\R;L^2((0,1);\R^2))$, that the elements of $W^\gamma$ satisfy (\ref{per}), 
that the elements of $V^\gamma(r_0,r_1)$ satisfy (\ref{per}) and (\ref{bc}) and such that
the left hand side of (\ref{eq:1.1}) is a Fredholm operator of index zero from $V^\gamma(r_0,r_1)$
into $W^\gamma$.

The main tools of the proofs are separation of variables 
(cf. (\ref{eq:3.3})--(\ref{eq:3.4})), integral
representation of the solutions of the corresponding boundary value problems 
of the ODE systems (cf. (\ref{eq:3.7}))
and an abstract criterion for Fredholmness which seems to be new (cf. 
Lemma \ref{lem:4.1}).

In order to formulate our results exactly, let us introduce function 
spaces: 
For $l=0,1,2,\dots$ and $\ga\ge 0$ we denote by $H^{l,\ga}$ the vector 
space of all
measurable functions $u: [0,1]\times\R\to\R$ such that
$
u(x,t)=u\left(x,t+\frac{2\pi}{\om}\right) \,\,\mbox{for\,\,almost\,\,all}\,
\,(x,t)\in[0,1]\times\R
$
and that
\begin{equation}\label{eq:1.3}
\|u\|_{H^{l,\ga}}^2:=\sum\limits_{k\in\Z}(1+k^2)^{\gamma}
\sum\limits_{m=0}^l\int\limits_0^1\left|\int\limits_0^{\frac{2\pi}{\om}}
\partial_x^mu(x,t)e^{-ik\om t}\,dt\right|^2\,dx<\infty.
\end{equation}
It is well-known (see, e.g., \cite{herrmann}, \cite[Chapter 
5.10]{robinson}, 
and \cite[Chapter 2.4]{vejvoda}) 
that $H^{l,\ga}$ is a Banach space. 
In fact, it is the space of all locally quadratically Bochner
integrable maps $u: \R \to H^l(0,1)$, which are $\frac{2 \pi}{\om}$--periodic 
and such that 
all generalized derivatives up to the (possibly noninteger) order $\ga$ 
are locally quadratically 
integrable maps from $\R$ into $H^l(0,1)$, too.
Further, we denote
\begin{eqnarray*}
W^{\gamma}&:=&H^{0,\ga}\times H^{0,\ga},\\
V^{\gamma}&:=& \left\{(u,v)\in W^{\gamma}:\,
(\d_tu+\d_xu,\d_tv-\d_xv)\in W^{\gamma}\right\}.
\end{eqnarray*}
The function spaces $W^{\ga}$ and $V^{\gamma}$ will be endowed with the norms
$$
\|(u,v)\|_{W^{\gamma}}^2:=\|u\|_{H^{0,\ga}}^2+\|v\|_{H^{0,\ga}}^2
$$
and 
$$
\|(u,v)\|_{V^{\gamma}}^2:=\|(u,v)\|_{W^{\gamma}}^2
+\|(\partial_tu+\partial_xu,\partial_tv-\partial_xv)\|_{W^{\gamma}}^2.
$$
It is easy to prove (see Lemma \ref{lem:2.2}) that $V^\ga$ is a Banach 
space. Moreover, 
if $\ga \ge 1$, then for any $(u,v) \in V^\ga$ and any $ x \in [0,1]$
there exist continuous trace maps $ (u,v) \in V^\ga \mapsto \left(u(x,\cdot), v(x,\cdot)\right) \in
\left(L_{loc}^2(\R)\right)^2$ (see Remark \ref{rem:2.4}). Hence, for $\ga \ge 
1$ it makes sense
to consider the following closed subspaces in $V^\ga$:
\begin{eqnarray*}
V^\ga(r_0,r_1)&:=&\{(u,v) \in V^\ga: \; (\ref{bc}) \mbox{ is fulfilled for a.a. } t \in \R\},\\
\tilde{V}^\ga(r_0,r_1)&:=&\{(u,v) \in V^\ga: \; (\ref{adbc}
) \mbox{ is fulfilled for a.a. } t \in \R\}.
\end{eqnarray*}
Finally, let us introduce linear operators: For $a,b,c,d\in L^{\infty}(0,
1)$ we define 
$A\in\LL(V^{\gamma}(r_0,r_1);W^{\gamma})$,
$\tilde{A}\in\LL(\tilde{V}^{\gamma}(r_0,r_1);W^{\gamma})$,
and $B,\tilde{B}\in\LL(W^{\gamma})$ by
$$
A
\left[
\begin{array}{c}
u\\v
\end{array}
\right]
:=
\left[
\begin{array}{c}
\partial_tu+\partial_xu+au\\\partial_tv-\partial_xv+dv
\end{array}
\right], \;\;\;
B
\left[
\begin{array}{c}
u\\v
\end{array}
\right]
:=
\left[
\begin{array}{c}
bv\\cu
\end{array}
\right],
$$
$$
\tilde{A}
\left[
\begin{array}{c}
u\\v
\end{array}
\right]
:=
\left[
\begin{array}{c}
-\partial_tu-\partial_xu+au\\
-\partial_tv+\partial_xv+dv
\end{array}
\right], \;\;\;
\tilde{B}
\left[
\begin{array}{c}
u\\v
\end{array}
\right]
:=
\left[
\begin{array}{c}
cv\\bu
\end{array}
\right].
$$
~\\

Now we formulate our main result:
\begin{thm}\label{thm:fredh}
Let $\ga \ge 1$, $a,d\in\L^{\infty}(0,1)$, $b,c\in BV(0,1)$, and suppose
(\ref{eq:1.4}).
Then we have:

(i) The operator $A$ is an isomorphism from $V^{\gamma}(r_0,r_1)$ onto 
$W^{\gamma}$.

(ii) The operator $A+B$ is Fredholm of index zero
from $V^{\gamma}(r_0,r_1)$ into $W^{\gamma}$.

(iii) The image of $A+B$ is the set of all $(f,g) \in W^\ga$ such that
$$
\int\limits_0^{\frac{2\pi}{\om}}
\int\limits_0^1\left(f(x,t)u(x,t)+g(x,t)v(x,t)\right)\,dx\,dt=0
\,\,\mbox{for \,all}\,\, (u,v)\in\ker(\tilde A+\tilde B).
$$
\end{thm}

In Theorem \ref{thm:fredh} and in what follows we denote, as usual, by 
$BV(0,1)$ 
the Banach space of all functions $h:(0,1) \to \R$ with bounded variation,
i.e. of
all $h \in L^\infty(0,1)$ such that there exists $C>0$ with
\beq
\label{BV}
\left|\int_0^1 h(x)\vphi'(x) dx \right| \le C \|\vphi\|_{L^\infty(0,
1)} \mbox{ for all } 
\vphi \in C_0^\infty(0,1).
\ee
The norm of $h$ in  $BV(0,1)$ is the sum of the norm of  $h$ in  $L^\infty(0,
1)$ and of the smallest 
possible constant $C$ in (\ref{BV}).\\

The present paper has been motivated mainly by two reasons:

The first reason are recent investigations of so-called 
traveling wave models, 
which are successfully used for modeling of dynamical behavior 
of semiconductor lasers (see, e.g., \cite{brs,LiRaRe,Peter,Rad,RadWu,Sieber}).
In those models,
systems of the type (\ref{eq:1.1}) appear to describe the forward and 
backward traveling light waves 
in the longitudinal laser direction. If one deals with multisection lasers 
consisting of several sections with different electrical and 
optical properties, then the coefficient functions $a$, $b$, $c$, and $d$ are 
discontinuous (piecewise constant).
Note that in the laser models the coefficient functions
and the reflection coefficients $r_0$ and $r_1$ as well as the 
unknown functions $u$ and $v$ are complex valued, and the linear wave 
system is coupled 
(via the coefficient functions) to a nonlinear balance equation for the 
carrier distribution 
in the active zone of the laser. Hence, from the point of view of applications 
to laser dynamics,  
our problem~(\ref{eq:1.1})--(\ref{bc}) is only a 
case study. In a forthcoming paper we will use 
our present results  for
applications to laser models.

The second reason is that the Fredholm property
of the linearization 
is  a key for many local investigation techniques for nonlinear equations, such as 
smooth continuation via implicit function theorem or local bifurcation
via Liapunov-Schmidt procedure. In particular, those techniques are well established 
for periodic solutions to 
nonlinear ODEs (see, e.g., \cite{IoJo}) and nonlinear parabolic PDEs
(see, e.g., \cite{Ki}). But almost nothing is known about the question if those techniques 
work for nonlinear dissipative hyperbolic PDEs.

\subsection{Some Remarks}\label{sec:remarks}

In this subsection we will comment about several aspects of Theorem \ref{thm:fredh}.

\begin{rem}\label{rem:maxreg} {\bf about maximal regularity} 
The main concern of this paper is not to prove existence of solutions to~(\ref{eq:1.1})--(\ref{bc}) 
for reasonable right hand sides $(f,g)$. The main concern is to find pairs of maximal regularity for 
(\ref{eq:1.1})--(\ref{bc}), i.e. Banach spaces $V$ and $W$ (in our case $V=V^\gamma(r_0,r_1)$
and $W=W^\gamma$) such that, on the one hand,  for all $(u,v) \in V$  (\ref{per}) and (\ref{bc}) are satisfied 
and the left hand side of (\ref{eq:1.1}) belongs to $W^\gamma$ and, on the other hand, 
for all $(f,g) \in W^\gamma$ the solutions to~(\ref{eq:1.1})--(\ref{bc}) belong to $V^\gamma(r_0,r_1)$.
\end{rem}

\begin{rem}\label{rem:choice} {\bf about the choice of the spaces} 
Our strategy of construction of  pairs of maximal regularity  for 
(\ref{eq:1.1})--(\ref{bc}) is as follows: 
First take a scale of spaces for the right hand 
sides, in our case $W^\gamma=H^{0,\gamma} \times H^{0,\gamma}$. 
Then take the maximal domain of definition of the differential
operator $(\partial_t+\partial_x,\partial_t-\partial_x)$ in $W^\gamma$, i.e. the space of all 
$(u,v) \in W^\gamma$ such that $(\partial_tu+\partial_xu,\partial_tv-\partial_xv) \in W^\gamma$.
Remark that this space is larger than the space of all $(u,v) \in W^\gamma$ such that
$\partial_tu,\partial_xu,\partial_tv,\partial_xv \in H^{0,\gamma}$.
In other words: The time derivatives and the space derivatives of $(u,v) \in V^\gamma$
can have singularities, such that  
$\partial_tu,\partial_xu,\partial_tv,\partial_xv \notin H^{0,\gamma}$,  
but these singularities cancel each other in $\d_tu+\d_xu$ and $\d_tv-\d_xv$.
In particular, it seems not to be a good idea to try to work with the 
time derivative operator and the space derivative operator separately.

And finally, take the scale parameter $\gamma$ large enough such that the elements of
the maximal domain of definition of the differential
operator $(\partial_t+\partial_x,\partial_t-\partial_x)$ have traces in $x=0,1$.
\end{rem}

\begin{rem}\label{rem:but} {\bf about noncompactness of the embedding 
\boldmath$V^\gamma(r_0,r_1) \hookrightarrow W^\gamma$\unboldmath}
At first glance it seems to be natural to prove assertion (ii) of Theorem \ref{thm:fredh} by using
 assertion (i) of Theorem \ref{thm:fredh} and by proving that $V^\gamma(r_0,r_1)$ is
compactly embedded into $W^\gamma$. But this approach fails because  $V^\gamma(r_0,r_1)$ is not
compactly embedded into $W^\gamma$! Let us verify this, for the sake of simplicity, for the case 
$r_0=r_1=0$: Take the sequence $(u_l,v_l) \in V^\gamma(0,0), l=1,2,\ldots$,
$$
u_l(x,t):=(1+l^2)^{-\frac{\gamma}{2}}e^{il\om(t-x)}x,\;
v_l(x,t):=(1+l^2)^{-\frac{\gamma}{2}}e^{il\om(t+x)}(x-1).
$$
Then
$$
\left(\partial_t+\partial_x\right)u_l(x,t)=(1+l^2)^{-\frac{\gamma}{2}}e^{il\om(t-x)},\;
\left(\partial_t-\partial_x\right)v_l(x,t)=(1+l^2)^{-\frac{\gamma}{2}}e^{il\om(t+x)}.
$$
Hence, $(u_l,v_l)$ is a bounded sequence in $W^\gamma$, but the sequence $(u_l,v_l)$
does not contain a subsequence which converges in $V^\gamma(0,0)$.
\end{rem}

\begin{rem}\label{rem:Aminus} {\bf about noncompactness of \boldmath$A^{-1}$\unboldmath and compactness
 of \boldmath$(A^{-1})^2$\unboldmath}
Remark \ref{rem:but} shows that the operator $A^{-1}$ is not compact from $W^\gamma$ into 
$W^\gamma$. But in Section 4 we will show that $(A^{-1})^2$ 
is compact from $W^\gamma$ into  $W^\gamma$. Using that and a corresponding abstract criterion for Fredholmness
(Lemma \ref{lem:4.1}), we will prove assertion  (ii) of Theorem \ref{thm:fredh}.
\end{rem}

\section{Some properties of the function spaces}\label{sec:spaces}
\renewcommand{\theequation}{{\thesection}.\arabic{equation}}
\setcounter{equation}{0}

In this section we formulate and prove some properties of the function spaces
$H^{l,\gamma}$, $W^{\ga}$, and $V^{\ga}$, introduced in Section~1.

It is well known that for each $u\in H^{l,\gamma}$ we have
\begin{equation}\label{eq:2.1}
u(x,t)=\sum\limits_{k\in\Z}u_k(x)e^{ik\omega t}\,\,\mbox{with}\,\,
u_k(x)=\frac{\om}{2\pi}\int\limits_0^{\frac{2\pi}{\om}}u(x,t)e^{-ik\omega 
t}\,dt,
\end{equation}
where the Fourier coefficients $u_k$ belong to the classical Sobolev space 
$H^l((0,1);\C)$, and the series
in~(\ref{eq:2.1}) converges in the complexification of
$H^{l,\gamma}$. And vice versa: For any sequence $(u_k)_{k\in\Z}$
with
\begin{equation}
\label{eq:2.2}
u_k\in H^l((0,1);\C),\;\;\overline{u_k}=u_{-k},\;\; 
\sum\limits_{k\in\Z}(1+k^2)^{\ga}\|u_k\|_{H^l((0,1);\C)}^2<\infty
\end{equation}
there exists exactly one $u\in H^{l,\ga}$ with~(\ref{eq:2.1}). In what 
follows, we will identify
functions $u\in H^{l,\gamma}$ and sequences $(u_k)_{k\in\Z}$ with~(\ref{eq:2.2})
by means of~(\ref{eq:2.1}), and we will keep 
for corresponding functions and sequences the notations $u$ and $(u_k)_{k\in\Z}$,
respectively.

\begin{lemma}\label{lem:2.1}
 A set $M\subset H^{0,\gamma}$ is precompact in $H^{0,\gamma}$ if and 
only if the 
following two conditions are satisfied:

(i) Uniform boundedness: There exists $C>0$ such that for all $u\in M$ it holds
 $$
\sum\limits_{k\in\Z}(1+k^2)^{\gamma}\int_0^1|u_k(x)|^2\,dx\le C.
$$

(ii) Uniform continuity with respect to shifts:
For all $\eps>0$ there exists $\de>0$ such that for all $\xi,\tau\in(-
\de,\de)$ and
all $u\in M$ it holds
$$
\sum\limits_{k\in\Z}(1+k^2)^{\gamma}\int\limits_{0}^1\Bigl|u_k(x+\xi)e^{ik\om 
\tau}-u_k(x)\Bigr|^2\,dx<
\eps,
$$
where $u_k(x+\xi):=0$ for $x+\xi\not\in[0,1]$.
\end{lemma}

\begin{proof}
Let us use the canonical isomorphism $J$ from $H^{0,\gamma}$ onto
$H^{0,0}=L^2\Bigl((0,1)\times\left(0,\frac{2\pi}{\om}\right)\Bigr)$, which 
is defined by 
$$
(Ju)(x,t):=\sum\limits_{k\in\Z}(1+k^2)^{\gamma/2}u_k(x)e^{ik\om t}.
$$
We have to show that $J(M)$ is precompact in $L^2\Bigl((0,1)\times\left(0,
\frac{2\pi}{\om}\right)\Bigr)$,
i.e. that $J(M)$ is bounded in $L^2\Bigl((0,1)\times\left(0,\frac{2\pi}{\om}\right)\Bigr)$,
and that
for all $\eps>0$ there exists $\de>0$ such that for all $\xi, \tau \in (-\delta,\delta)$ and all $u\in J(M)$ it holds
\begin{equation}\label{eq:2.3}
\int\limits_{0}^{\frac{2\pi}{\om}}\int\limits_{0}^1|u(x+\xi,
t+\tau)-u(x,t)|^2\,dx\,dt<
\eps,
\end{equation}
where $u(x+\xi,t+\tau):=0$ for $x+\xi\not\in[0,1]$. Boundedness in 
$L^2\Bigl((0,1)\times\left(0,\frac{2\pi}{\om}\right)\Bigr)$ is just condition 
$(i)$ of the lemma.

Now we show that~(\ref{eq:2.3})  is just condition $(ii)$ of the lemma.
This follows from
\begin{eqnarray*}
\lefteqn{
\int\limits_{0}^{\frac{2\pi}{\om}}\int\limits_{0}^1\Bigl|(Ju)(x+\xi,
t+\tau)-(Ju)(x,t)\Bigr|^2\,dx\,dt}\\
&&=\int\limits_{0}^{\frac{2\pi}{\om}}\int\limits_{0}
^1\left|\sum\limits_{k\in\Z}(1+k^2)^{\gamma/2}
\left(u_k(x+\xi)e^{ik\om(t+\tau)}-u_k(x)e^{ik\om t}\right)\right|^2\,dx\,dt\\
&&=\frac{2\pi}{\om}\int\limits_{0}^1\sum\limits_{k\in\Z}(1+k^2)^{\gamma}
\Bigl|u_k(x+\xi)e^{ik\om \tau}-
u_k(x)\Bigr|^2\,dx.
\end{eqnarray*}
\end{proof}

\begin{lemma}\label{lem:2.2}
The space $V^{\gamma}$ is complete.
\end{lemma}

\begin{proof}
Let $(u^{j},v^{j})_{j\in\N}$ be a fundamental sequence 
in $V^{\gamma}$. 
Then $(u^j,v^j)_{j\in\N}$ and $(\d_tu^j+\d_xu^j,\d_tv^j-\d_xv^j)_{j\in\N}$
are fundamental sequences in $W^{\gamma}$. Because $W^{\gamma}$ is complete,
there exist
$(u,v)\in W^{\gamma}$ and $(\tilde u,\tilde v)\in W^{\gamma}$ such that
$$
(u^{j},v^{j})\to (u,v) \; \mbox{ and } \;
(\d_tu^j+\d_xu^j,\d_tv^j-\d_xv^j)\to (\tilde u,\tilde v)
$$
in $W^{\ga}$ as $j\to\infty$.
It remains to show that
$
\d_tu+\d_xu=\tilde u
$
and
$
\d_tv-\d_xv=\tilde v
$
in the sense of generalized derivatives. But this is obvious: Take a smooth 
function 
$\vphi: (0,1)\times\left(0,\frac{2\pi}{\om}\right)\to\R$ with compact 
support. Then
\begin{eqnarray*}
&\displaystyle
\int\limits_{0}^{\frac{2\pi}{\om}}\int\limits_{0}^1u(\d_t+\d_x)\vphi\,
dx\,dt=
\lim\limits_{j\to\infty}\int\limits_{0}^{\frac{2\pi}{\om}}\int\limits_{
0}^1
u^j(\d_t+\d_x)\vphi\,dx\,dt&\\
&\displaystyle
=-\lim\limits_{j\to\infty}\int\limits_{0}^{\frac{2\pi}{\om}}\int\limits_{
0}^1
(\d_t+\d_x)u^j\vphi\,dx\,dt=-
\int\limits_{0}^{\frac{2\pi}{\om}}\int\limits_{0}^1
\tilde u\vphi\,dx\,dt,&
\end{eqnarray*}
and similarly for $v$ and $\tilde v$.
\end{proof}

\begin{lemma}\label{lem:2.3}
If $\gamma\ge 1$, then  $V^{\gamma}$ is continuously embedded into $\left(H^{1,
\ga-1}\right)^2$.
\end{lemma}

\begin{proof}
Take $(u,v)\in V^{\ga}$. Then $(u,v)\in\left(H^{0,\ga}\right)^2$, 
hence $(\d_tu,\d_tv)\in\left(H^{0,\ga-1}\right)^2$.
By the definition of the space $V^{\ga}$, $(\d_xu,\d_xv)\in\left(H^{0,
\ga-1}\right)^2$.
Hence $(u,v)\in\left(H^{1,\ga-1}\right)^2$. Moreover, we have
\begin{eqnarray*}
\lefteqn{
\|(u,v)\|_{\left(H^{1,\ga-1}\right)^2}^2=\|u\|_{H^{0,\ga-1}}^2+\|v\|_{H^{0,
\ga-1}}^2+
\|\d_xu\|_{H^{0,\ga-1}}^2+\|\d_xv\|_{H^{0,\ga-1}}^2}\\
&&\le 
\|u\|_{H^{0,\ga-1}}^2+\|v\|_{H^{0,\ga-1}}^2+
\|\d_tu+\d_xu\|_{H^{0,\ga-1}}^2+\|\d_tv-\d_xv\|_{H^{0,\ga-1}}^2\\
&&\;\;\;\;\;\;\;\;\;\;\;\;\;\;\;\;+\|\d_tu\|_{H^{0,\ga-1}}^2+\|\d_tv\|_{H^{0,
\ga-1}}^2 \\
&&\le C\|(u,v)\|^2_{V^{\gamma}},
\end{eqnarray*}
where the constant $C$ does not depend on $(u,v)$.
\end{proof}

\begin{rem}\label{rem:2.4}
Suppose $\ga\ge 1$. By Lemma \ref{lem:2.3}, we have
$$
V^\ga \hookrightarrow \left(H^{1,0}\right)^2 
\approx \left(H^1\left((0,1);L^2\left(0,\frac{2\pi}{\om}\right)\right)\right)^2
\hookrightarrow \left(C\left([0,1];L^2\left(0,\frac{2\pi}{\om}\right)\right)\right)^2.
$$
Therefore, for any $x \in [0,1]$ there exists a continuous trace map
$$
(u,v) \in V^\ga \mapsto \Bigl(u(x,\cdot),v(x,\cdot)\Bigr) 
\in \left(L^2\left(0,\frac{2\pi}{\om}\right)\right)^2.
$$
\end{rem}

Now, let us consider the dual spaces $(H^{0,\ga})^*$.

Obviously, for any $\ga \ge 0$ the spaces $H^{0,\ga}$ are densely and 
continuously embedded into
the Hilbert space $H^{0,0}=L^2\left((0,1) \times \left(0,\frac{2 \pi}{\om}\right)\right)$.
Hence, there is a 
canonical dense continuous embedding
\beq
\label{emb}
H^{0,0} \hookrightarrow (H^{0,\ga})^*: \;\;[u,v]_{H^{0,\ga}} = \langle 
u,v \rangle_{H^{0,0}}
\mbox{ for all } u \in H^{0,0}  \mbox{ and }  v \in H^{0,\ga}.
\ee
Here $[\cdot,\cdot]_{H^{0,\ga}}:(H^{0,\ga})^* \times H^{0,\ga} \to \R$ 
is the dual pairing,
and $\langle \cdot,\cdot \rangle_{H^{0,0}}: H^{0,0} \times H^{0,0} \to 
\R$ is the scalar product in 
$H^{0,0}$, i.e.
\beq
\label{SP}
\langle u,v \rangle_{H^{0,0}}:=\frac{\om}{2 \pi}
\int\limits_{0}^\frac{2\pi}{\om}\int\limits_{0}^{
1}u(x,t)v(x,t)\,dx\,dt=
\sum_{k \in \Z} \int\limits_{0}^1 u_k(x) \overline{v_k(x)}\,dx.
\ee
Let us denote
\beq
\label{e}
e_k(t):=e^{ik\om t} \mbox{ for } k \in \Z \mbox{ and } t \in \R.
\ee
If a sequence $(\vphi_k)_{k \in \Z}$ with $\vphi_k \in L^2((0,1);\C)$ 
is given, then the pointwise
products $\vphi_ke_k$ belong to $L^2\left((0,1)\times \left(0,\frac{2\pi}{\om}\right);
\C\right)$. 
Hence, they belong to the complexification of 
$(H^{0,\ga})^*$ (by means of the complexification of (\ref{emb})), and it makes sense to ask if the series
\beq
\label{series}
\sum_{k \in \Z} \vphi_ke_k
\ee
converges in the complexification of $(H^{0,\ga})^*$.

\begin{lemma}\label{lem:2.5}
(i) \, For any $\vphi \in (H^{0,\ga})^*$ there exists a sequence $(\vphi_k)_{k 
\in \Z}$ with
\beq
\label{seqBed}
\vphi_k \in L^2((0,1);\C), \; \overline{\vphi_k}=\vphi_{-k}, \; 
\sum_{k \in \Z} (1+k^2)^{-\ga}\int\limits_{0}^1 |\vphi_k(x)|^2 
\,dx < \infty,
\ee
such that  the series (\ref{series}) converges to $\vphi$ in 
the complexification of $(H^{0,\ga})^*$.
Moreover, it holds
\beq
\label{coeffBed}
\int\limits_{0}^1 \vphi_k(x)u(x) \,dx =[\vphi,ue_{-k}]_{H^{0,
\ga}} \mbox{ for all  } u \in L^2(0,1).
\ee

(ii) \, For any sequence $(\vphi_k)_{k \in \Z}$ with (\ref{seqBed}) the series (\ref{series}) converges
in the complexification of 
$(H^{0,\ga})^*$ to some $\vphi \in (H^{0,\ga})^*$, and (\ref{coeffBed}) is satisfied.
\end{lemma}

\begin{proof}(i) \, By the Riesz representation theorem, for 
given $\vphi \in (H^{0,\ga})^*$
and $k \in \Z$ there exists exactly one $\vphi_k \in L^2((0,1);\C)$ with 
(\ref{coeffBed}). 
The property $\overline{\vphi_k}=\vphi_{-k}$ follows directly from  (\ref{coeffBed}).

Now, take $u \in H^{0,\ga}$ with its representation (\ref{eq:2.1}), (\ref{eq:2.2}) 
(with $l=0$ there). We have
$$
\Bigl|[\vphi,u]_{H^{0,\ga}}\Bigr|=\left|\left[\vphi,\sum_{k \in \Z} u_ke_k\right]_{H^{0,
\ga}}\right|=
\left|\sum_{k \in \Z}\int\limits_{0}^1 (1+k^2)^{-\ga}\vphi_k(x) 
(1+k^2)^\ga 
\overline{u_k(x)} \,dx \right|.
$$
Taking the supremum over all $u \in H^{0,\ga}$ with
$
\|u\|_{H^{0,\ga}}^2=\sum_{k \in \Z} (1+k^2)^\ga \int_{0}^1 |u_k(x)|^2 
\,dx =1,
$
we get
$$
\|\vphi\|^2_{(H^{0,\ga})^*}=\sum_{k \in \Z} (1+k^2)^{-\ga} \int\limits_{
0}^1 |\vphi_k(x)|^2 \,dx < \infty.
$$
Similarly one shows that the series (\ref{series}) converges 
in the complexification of $(H^{0,\ga})^*$:
For all $u \in H^{0,\ga}$ we have
\begin{eqnarray}
\label{absch}
\lefteqn{
\left|\left[\vphi-\sum_{|k| \le k_0} \vphi_ke_k,u\right]_{H^{0,\ga}}\right|
=\left|\left[\vphi-\sum_{|k| \le k_0} \vphi_ke_k,\sum_{k \in \Z} u_ke_k\right]_{H^{0,
\ga}}\right|}
\nonumber\\
&&=\left|\sum_{|k| > k_0}\int\limits_{0}^1 (1+k^2)^{-\ga}\vphi_k(x) 
(1+k^2)^\ga 
\overline{u_k(x)} \,dx \right|
\nonumber\\
&&\le \left(\sum_{|k| > k_0}(1+k^2)^{-\ga}\int\limits_{0}^1 
|\vphi_k(x)|^2 \; dx \right)^{\frac{1}{2}}
\|u\|_{H^{0,\ga}},
\end{eqnarray}
and this tends to zero as $k_0 \to \infty$ uniformly for $\|u\|_{H^{0,\ga}}=1$.

(ii) \, As in (\ref{absch}) one shows that (\ref{series})
 converges in the complexification 
of $(H^{0,\ga})^*$
to some $\vphi \in (H^{0,\ga})^*$. Moreover, (\ref{coeffBed})
 is satisfied, because
for all $u \in L^2(0,1)$ we have
$$
[\vphi,ue_{-k}]_{H^{0,\ga}}=\left[\sum_{l \in \Z}\vphi_l e_l,ue_{-k}\right]_{H^{0,
\ga}}
=\sum_{l \in \Z}\langle \vphi_l e_l,ue_k \rangle_{H^{0,0}}
=\int\limits_{0}^1 \vphi_k(x) u(x) \,dx.
$$
Here we used (\ref{emb}) and (\ref{SP}).
\end{proof}

\section{Proof of the isomorphism property}\label{sec:iso}
\setcounter{equation}{0}

In what follows, we suppose the assumptions of 
Theorem~\ref{thm:fredh} to be fulfilled. 
In this section we prove assertion $(i)$ of Theorem~\ref{thm:fredh}. 

Fix $(f,g)\in W^{\ga}$. Then
$f(x,t)=\sum\limits_{k\in\Z}f_k(x)e^{ik\omega t}$ and 
$g(x,t)=\sum\limits_{k\in\Z}g_k(x)e^{ik\omega t}$
with $f_k,g_k\in\L^2\left((0,1);\C\right)$ and
\begin{equation}\label{eq:3.1}
\sum\limits_{k\in\Z}(1+k^2)^{\gamma}\int\limits_{0}^1 \left(|f_k(x)|^2+|g_k(x)|^2\right)\,
dx<\infty.
\end{equation}
We have to show that there exists exactly one $(u,v)\in V^{\ga}(r_0,r_1)$ 
such that
$\d_tu+\d_xu+au=f$ and $\d_tv-\d_xv+du=g.$
Writing $u$ and $v$ as series according to~(\ref{eq:2.1}) and~(\ref{eq:2.2}),  
we have to show 
that there exists exactly one pair of sequences $(u_k)_{k\in\Z}$ and
$(v_k)_{k\in\Z}$ with $u_k,v_k\in H^1(0,1)$ satisfying the boundary value 
problem
\begin{equation}\label{eq:3.3}
u_k^{\prime}+(a(x)+ik\omega)u_k=f_k(x),\;
v_k^{\prime}-(d(x)+ik\omega)v_k=-g_k(x),
\end{equation}
\begin{equation}\label{eq:3.4}
u_k(0)=r_0v_k(0),\;
v_k(1)=r_1u_k(1),
\end{equation}
and the estimates
\begin{equation}\label{eq:3.5}
\sum\limits_{k\in\Z}(1+k^2)^{\gamma}\int\limits_{0}^1\left(|u_k(x)|^2+
|v_k(x)|^2\right)dx <\infty,
\end{equation}
\begin{equation}\label{eq:3.6}
\sum\limits_{k\in\Z}(1+k^2)^{\gamma}\int\limits_{0}^1\left(|u_k^{\prime}(x)+ik\omega 
u_k(x)|^2+
|v_k^{\prime}(x)-ik\omega v_k(x)|^2\right) dx <\infty.
\end{equation}
Here we used Lemma~\ref{lem:2.3} and Remark~\ref{rem:2.4}.

The estimate~(\ref{eq:3.6}) follows from~(\ref{eq:3.1}),~(\ref{eq:3.3}),
and~(\ref{eq:3.5}).
Hence, it remains to show that there exists exactly one pair of sequences
$(u_k)_{k\in\Z}$ and $(v_k)_{k\in\Z}$ with $u_k,v_k\in H^1\left((0,1);
\C\right)$ 
satisfying~(\ref{eq:3.3}),~(\ref{eq:3.4}), and~(\ref{eq:3.5}).

In order to simplify the formulae below, let us introduce the following 
notation:
$$
\al(x):=\int\limits_0^xa(y)\,dy,\quad \de(x):=\int\limits_0^xd(y)\,dy,\quad
\Delta_k:=e^{ik\om+\de(1)}-r_0r_1e^{-ik\om-\al(1)}.
$$
A straightforward calculation shows that the boundary value problem~(\ref{eq:3.3}),~(\ref{eq:3.4})
has a unique solution $(u_k,v_k)\in H^1\left((0,1);\C^2\right)$, 
and this solution is explicitely given by
\begin{equation}\label{eq:3.7}
\begin{array}{rcl}
\displaystyle
u_k(x)&=&e^{-ik\om x-\al(x)}\left(
\displaystyle\int\limits_0^xe^{ik\om y+\al(y)}f_k(y)\,dy
+\frac{r_0}{\Delta_k}w_k(f_k,g_k)\right),\\
v_k(x)&=&e^{ik\om x+\de(x)}\left(
\displaystyle\int\limits_0^xe^{-ik\om y-\de(y)}g_k(y)\,dy
+\frac{1}{\Delta_k}w_k(f_k,g_k)\right)
\end{array}
\end{equation}
with
\begin{equation}\label{eq:w_k}
w_k(f,g):=r_1e^{-ik\om-\al(1)}
\int\limits_0^1e^{ik\om y+\al(y)}f(y)\,dy-
e^{ik\om+\de(1)}
\int\limits_0^1e^{-ik\om y-\de(y)}g(y)\,dy.
\end{equation}
Here we used assumption~(\ref{eq:1.4}), which implies
\begin{equation}\label{eq:3.9}
|\Delta_k|\ge\Bigl|e^{\de(1)}-|r_0r_1|e^{-\al(1)}\Bigr|>0
\mbox{ for all } k \in \Z.
\end{equation}
From~(\ref{eq:3.7}) and~(\ref{eq:3.9}) it follows that
\begin{equation}\label{eq:3.10}
|u_k(x)|+|v_k(x)|\le C\left(\int\limits_0^1\left(|f_k(x)|^2+|g_k(x)|^2 
\right) dx \right)^\frac{1}{2}
\end{equation}
for all $x\in [0,1]$, where the constant $C$ does not depend on $k$, $f_k$,
$g_k$, and $x$.
Finally,~(\ref{eq:3.1}) and~(\ref{eq:3.10}) imply~(\ref{eq:3.5}).

\section{Proof of Fredholmness}\label{sec:ferdh}
\setcounter{equation}{0}

In this section we prove that $A+B$ is Fredholm, which is a part of assertion
$(ii)$ of Theorem~\ref{thm:fredh}. 

Obviously, $A+B$ is Fredholm from $V^{\ga}(r_0,r_1)$ into $W^{\ga}$ if and
only if $I+BA^{-1}$ is Fredholm from $W^{\ga}$ into $W^{\ga}$. Here
$I$ is the identity in~$W^{\ga}$.

We are going to prove that $I+BA^{-1}$ is Fredholm from $W^{\ga}$ into 
$W^{\ga}$
using the following 

\begin{lemma}\label{lem:4.1}
Let $W$ be a Banach space, $I$ the identity in $W$, and $C\in\LL(W)$ such 
that $C^2$
is compact.  Then $I+C$ is Fredholm.
\end{lemma}

\begin{proof}
Since $I-C^2=(I-C)(I+C)$ and $C^2$ is compact, we have
\beq
\label{endl}
\dim\ker(I+C)\le \dim\ker(I-C^2)<\infty.
\ee
Similarly one gets $\dim\ker(I+C)^*<\infty$, hence
$
\codim\overline{\im(I+C)}<\infty.
$
It remains to show that $\im(I+C)$ is closed.

Take a sequence $(w_j)_{j\in\N}\subset W$ and an element $w\in W$ such that
\begin{equation}\label{eq:4.1}
(I+C)w_j\to w.
\end{equation}
We have to show that $w \in \im(I+C)$.

Because of (\ref{endl}) there exists a closed subspace $V$ of 
$W$ such that 
\beq
\label{dirSum} 
W=\ker(I+C)\oplus V,
\ee
Using the decomposition
$w_{j}=u_j+v_j$ with  $u_j\in\ker(I+C)$ and $v_j \in V$,
we get from (\ref{eq:4.1})
\begin{equation}
\label{eq:v}
(I+C)v_j\to w.
\end{equation}

First we show that the sequence $(v_j)_{j \in \N}$ is bounded.
If not, without loss of generality we can assume that
\begin{equation}\label{eq:4.11}
\lim\limits_{j\to \infty}\|v_j\|=\infty.
\end{equation}
>From (\ref{eq:v}) and (\ref{eq:4.11}) we get 
\beq
\label{1Null}
(I+C)\frac{v_j}{\|v_j\|} \to 0,
\ee
hence
\beq
\label{Null}
(I-C^2)\frac{v_j}{\|v_j\|}  \to 0.
\ee
On the other side, because $C^2$ is compact, there exist $v \in W$ and 
a subsequence 
$(v_{j_k})_{k \in \N}$ such that
\beq
\label{conv}
C^2\frac{v_{j_k}}{\|v_{j_k}\|}  \to v.
\ee
Inserting (\ref{conv}) into (\ref{Null}), we get 
\beq
\label{minus}
\frac{v_{j_k}}{\|v_{j_k}\|} \to v \in V.
\ee
Combining (\ref{minus}) with (\ref{1Null}), we get 
$(I+C)v=0$,
i.e. 
$$v \in V \cap \ker (I+C) \mbox{ and } \|v\|=1.
$$ 
But this contradicts to (\ref{dirSum}).

Now we use the boundedness of $(v_j)_{j \in \N}$ to show that
that $w \in \im(I+C)$. 
As the operator $C^2$ is compact, there exist $v\in W$ and a subsequence 
$(v_{j_k})_{k \in \N}$ such that
$
C^2v_{j_k}\to v.
$
On the other hand, (\ref{eq:v}) yields $(I-C^2)v_j \to (I-C)w$. 
Hence
(\ref{eq:v}) yields
$
\lim\limits_{k\to \infty}v_{j_k}=(I-C)w+v,
$
and, therefore,
$$
w=\lim\limits_{k\to \infty}(I+C)v_{j_k}=(I+C)((I-C)w+v)\in\im(I+C).
$$
\end{proof}

In order to use Lemma~\ref{lem:4.1} with $W:=W^{\ga}$ and $C:=BA^{-1}$,
let us show that
$\left(BA^{-1}\right)^2$ is compact from $W^{\ga}$ into $W^{\ga}$.

Take a bounded set $N\subset W^{\ga}$, and let $M$ be its image under
$(BA^{-1})^2$. In order to show that $M$ is precompact in 
$W^{\ga}=H^{0,\ga}\times H^{0,\ga}$, we use Lemma~\ref{lem:2.1} ``componentwise''.

Condition $(i)$ of Lemma~\ref{lem:2.1} is satisfied because $(BA^{-1})^2$
is a bounded operator from $W^{\ga}$ into $W^{\ga}$.

It remains to check condition $(ii)$ of Lemma~\ref{lem:2.1}.
The explicit representation~(\ref{eq:3.7}) of $A^{-1}$ yields the following:
For given $(f,g)\in N$ we have
$$
\left[
\begin{array}{c}
\tilde u\\\tilde v
\end{array}
\right]
=\left(BA^{-1}\right)^2
\left[
\begin{array}{c}
f\\g
\end{array}
\right]
=BA^{-1}B
\left[
\begin{array}{c}
u\\v
\end{array}
\right]
\mbox{ with } 
\left[
\begin{array}{c}
u\\v
\end{array}
\right]=A^{-1}
\left[
\begin{array}{c}
f\\g
\end{array}
\right]
$$
if and only if 
$$
\tilde{u}_k(x)=b(x)e^{ik\om x+\de(x)}\left(
\int\limits_0^xe^{-ik\om y-\de(y)}c(y)u_k(y)\,dy
+\frac{1}{\Delta_k}w_k(bv_k,cu_k)\right),
$$
$$
\tilde{v}_k(x)=c(x)e^{-ik\om x-\al(x)}\left(
\int\limits_0^xe^{ik\om y+\al(y)}b(y)v_k(y)\,dy
+\frac{r_0}{\Delta_k}w_k(bv_k,cu_k)\right),
$$
where the functions $w_k$ are defined by~(\ref{eq:w_k}) and the functions 
$u_k,v_k\in H^1(0,1)$ 
as the solutions to~(\ref{eq:3.3}),~(\ref{eq:3.4})
are given by the formulas~(\ref{eq:3.7}). Hence
\begin{eqnarray*}
\lefteqn{\left|\tilde u_k(x+\xi)e^{ik\om \tau}-\tilde u_k(x)\right|}\\ 
&&\le\left|b(x+\xi)
e^{ik\om(\xi+\tau)+\de(x+\xi)}
\int\limits_x^{x+\xi}e^{-ik\om y-\de(y)}c(y)u_k(y)\,dy\right|\\
&&+\left|b(x+\xi)e^{\de(x+\xi)}\left(e^{2ik\om(\xi+\tau)}-1\right)\left(\int\limits_0^{
x}e^{-ik\om y-\de(y)}
c(y)u_k(y)\,dy 
+\frac{w_k(bv_k,cu_k)}{\Delta_k}\right)\right|\\
&&+\left|\left(b(x+\xi)e^{\de(x+\xi)}-b(x)e^{\de(x)}\right)\left(\int\limits_0^{
x}e^{-ik\om y-\de(y)}c(y)u_k(y)\,dy 
+\frac{w_k(bv_k,cu_k)}{\Delta_k}\right)\right|.
\end{eqnarray*}
On the account of~(\ref{eq:3.10}) and the boundedness of  $N$, we get 
the estimate
$$
\sum\limits_{k\in\Z}(1+k^2)^{\gamma}
\int\limits_{0}^1\left|b(x+\xi)
e^{ik\om(\xi+\tau)+\de(x+\xi)}
\int\limits_x^{x+\xi}e^{-ik\om y-\de(y)}c(y)u_k(y)\,dy\right|^2dx \le C\xi^2,
$$
where the constant $C$ does not depend on $\xi$, $\tau$, and $(f,g)\in N$.
Similarly one gets (using (\ref{eq:w_k})--(\ref{eq:3.10}))
\begin{eqnarray*}
\lefteqn{
\sum\limits_{k\in\Z}(1+k^2)^{\gamma}
\int\limits_{0}^1\left|b(x+\xi)e^{\de(x+\xi)}-b(x)e^{\de(x)}\right|^2}\\
&&\times\left|\int\limits_0^{x}e^{-ik\om y-\de(y)}c(y)u_k(y)\,dy 
+\frac{w_k(bv_k,cu_k)}{\Delta_k}\right|^2\,dx\\
&&\le C\int\limits_{0}^1\left|b(x+\xi)e^{\de(x+\xi)}-b(x)e^{\de(x)}\right|^2\,
dx,
\end{eqnarray*}
which tends to zero for $\xi \to 0$ uniformly with respect to $(f,g) \in N$
(because of the continuity in the mean of the function $x \mapsto b(x)e^{\de(x)}$).

It remains to show that
\begin{eqnarray}
\label{absch1}
\lefteqn{
\sum\limits_{k\in\Z}(1+k^2)^{\gamma}
\int\limits_{0}^1\left|
b(x+\xi)e^{\de(x+\xi)}\left(e^{ik\om(\xi+\tau)}-1\right)\right|^2}\nonumber\\
&&\times\left|\int\limits_0^{x}e^{-ik\om y-\de(y)}c(y)u_k(y)\,dy 
+\frac{w_k(bv_k,cu_k)}{\Delta_k}\right|^2dx
\end{eqnarray}
tends to zero for $\xi,\tau \to 0$ uniformly with respect to $(f,g) \in N$.
Using (\ref{eq:3.3}), we have for
all $k\ne 0$ 
\begin{eqnarray*}
\lefteqn{
e^{-ik\om y}c(y)u_k(y)}\\
&&
=\frac{1}{2ik\om}c(y)\left(e^{-2ik\om y}\frac{d}{
dy}
\left(e^{ik\om y}u_k(y)\right)-\frac{d}{dy}
\left(e^{-ik\om y}u_k(y)\right)\right)\\
&&
=\frac{1}{2ik\om}c(y)\left(e^{-ik\om y}
\left(f_k(y)-a(y)u_k(y)\right)-\frac{d}{dy}
\left(e^{ik\om y}u_k(y)\right)\right).
\end{eqnarray*}
Moreover, assumption $c\in BV(0,1)$ yields (cf. (\ref{BV}))
$$
\left|\int\limits_0^{x}e^{-\de(y)}c(y)\frac{d}{
dy}
\left(e^{ik\om y}u_k(y)\right)\,dy \right|\le C\|u_k\|_{L^{\infty}(0,1)},
$$
the constant $C$ being independent of $x$, $k$, and $u_k$. Using~(\ref{eq:3.10}), it follows
$$
\left|\int\limits_0^{x}e^{-ik\om y-\de(y)}c(y)u_k(y)\,dy \right|\le
\frac{C}{1+|k|} \left(\|f_k\|_{L^{2}(0,1)}+\|g_k\|_{L^{
2}(0,1)}\right)
$$
for some $C>0$ not depending
on $x$, $\xi$, $\tau$, $k$, $f_k$, and $g_k$. Similarly we proceed in
the integrals in $w_k(bv_k,cu_k)$ (cf. (\ref{eq:w_k})) in order to get
$$
\left|w_k(bv_k,cu_k)\right|\le
\frac{C}{1+|k|} \left(\|f_k\|_{L^{2}(0,1)}+\|g_k\|_{L^{
2}(0,1)}\right).
$$
Hence, (\ref{absch1}) can be estimated by 
$$
C\sum\limits_{k\in\Z}(1+k^2)^{\gamma-1}
\left|e^{ik\om(\xi+\tau)}-1\right|^2 \left(\|f_k\|^2_{L^{2}(0,1)}
+\|g_k\|^2_{L^{2}(0,1)}\right)
$$
with some $C>0$ not depending
on $x$, $\xi$, $\tau$, $k$, $f_k$, and $g_k$.
Using $|e^{ik\om(\xi+\tau)}-1|\le k\om(\xi+\tau)$, we see that
this tends to zero as $\xi,\tau \to 0$ uniformly with respect to $(f,g)\in N$.

Thus, $M$ is precompact, i.e. $(BA^{-1})^2$ is compact. Hence 
$I+BA^{-1}$ is Fredholm, and therefore $A+B$ is Fredholm.

In order to finish the proof of assertion $(ii)$ of Theorem~\ref{thm:fredh}, 
it remains to show that the index 
of $A+B$ is zero. This will be proved in the next section.

\section{Fredholmness of index zero}\label{sec:index}
\setcounter{equation}{0}


Directly from the definitions of the operators $A$, $\tilde{A}
$, $B$, and $\tilde{B}$ it follows
\begin{eqnarray}
\label{5.1}
\lefteqn{
\langle (A+B)(u,v),(\tilde{u},\tilde{v}) \rangle
=\langle (u,v),(\tilde{A}+\tilde{B})(\tilde{u}
,\tilde{v}) \rangle} \nonumber \\
&&\mbox{ for all } (u,v) \in V^\ga(r_0,r_1) \mbox{ and } (\tilde{
u},\tilde{v}) \in \tilde{V}^\ga(r_0,r_1).
\end{eqnarray}
Here 
$$
\langle (u,v),(\vphi,\psi)\rangle:=\frac{\om}{2 \pi} \int\limits_0^{\frac{2 
\pi}{\om}}
\int\limits_0^1\left(u\vphi+v\psi\right)\,dx\,dt =\sum\limits_{k \in \Z}\int\limits_0^1
\left(u_k\overline{\vphi}_k+v_k \overline{\psi}_k\right)\;dx
$$
is the usual scalar product in the Hilbert space $\left(H^{0,0}\right)^2=L^2\left((0,
1)\times
\left(0,\frac{2 \pi}{\om}\right);\R^2\right)$.

In order to prove that the Fredholm operator $A+B$ has index zero, it 
suffices to show that
\beq
\label{5.2}
\dim \ker (A+B)=\dim \ker (A+B)^*.
\ee
Here $(A+B)^*$ is the dual operator  to $A+B$, i.e. a linear bounded 
operator from $(W^\ga)^*$ into $(V^\ga(r_0,r_1))^*$. Using the continuous 
dense embedding
$$
\tilde{V}^\ga(r_0,r_1) \hookrightarrow W^\ga  \hookrightarrow 
(H^{0,0})^2  
\hookrightarrow \left(W^\ga\right)^*,
$$
it makes sense to compare the subspaces $\ker (A+B)^*$ of $\left(W^\ga\right)^*$ 
and $\ker (\tilde{A}+\tilde{B})$ of $\tilde{V}^\ga(r_0,r_1)$:

\begin{lemma}\label{lem:adj}
$\ker (A+B)^*=\ker (\tilde{A}+\tilde{B}).$
\end{lemma}

\begin{proof}
Let $[\cdot,\cdot]:(W^\ga)^* \times W^\ga \to \R$ be the dual pairing 
on $W^\ga$. Then for all 
$(u,v) \in V^\ga(r_0,r_1)$ and $(\tilde{u},\tilde{v}
) \in \tilde{V}^\ga(r_0,r_1)$ we have
\begin{eqnarray}
\label{5.3}
\lefteqn{\langle (\tilde{A}+\tilde{B})(\tilde{
u},\tilde{v}),(u,v)\rangle = 
\langle (\tilde{u},\tilde{v}),(A+B)(u,v)\rangle } 
\nonumber\\
&&=\left[(\tilde{u},\tilde{v}),(A+B)(u,v)\right]=
\left[(A+B)^*(\tilde{u},\tilde{v}),(u,v)\right].
\end{eqnarray}
Here we used (\ref{emb}) and (\ref{5.1}). Obviously, (\ref{5.3}) implies
$
\ker(\tilde{A}+\tilde{B})\subseteq \ker (A+B)^*.
$

Now, take an arbitrary $(\vphi,\psi) \in \ker(A+B)^*$, and let us show that
$(\vphi,\psi) \in \ker(\tilde{A}+\tilde{B})$.
By  Lemma \ref{lem:2.5}, we have (using notation (\ref{e}))
$
\vphi=\sum\limits_{k \in \Z} \vphi_ke_k$ and  $\psi=\sum\limits_{k \in 
\Z} \psi_ke_k
$
with
$
\vphi_k,\psi_k \in L^2\left((0,1);\C\right)
$
and
$$ 
\sum\limits_{k \in \Z}(1+k^2)^\ga
\int_0^1\left(|\vphi_k(x)|^2+ |\psi_k(x)|^2\right)dx < \infty.
$$ 
It follows  that for all $(u,v) \in V^\ga(r_0,r_1)$
\begin{eqnarray*}
\lefteqn{0=\left[(A+B)^*(\vphi,\psi),(u,v)\right]=
\left[(\vphi,\psi),(A+B)(u,v)\right]}\\
&&=\sum\limits_{k \in \Z} \int_0^1 (\vphi_k(u'_{-k}+
(a-ik\om)u_{-k}+bv_{-k})
+\psi_k(v'_{-k}-(d-ik\om)v_{-k}+cu_{-k}))dx.
\end{eqnarray*}
Therefore
$$
\vphi_k\left(u'_{-k}+(a-ik\om)u_{-k}+bv_{-k}\right)+
\psi_k\left(v'_{-k}-(d-ik\om)v_{-k}+cu_{-k}\right)
=0
$$
for all $u_k,v_k \in H^1\left((0,1);\C\right)$ with (\ref{eq:3.4}). By 
a standard argument, 
we conclude that
$\vphi_k,\psi_k \in H^1\left((0,1);\C\right)$ and that they satisfy the 
differential equations
\beq
\label{adjGl}
-\vphi'_k+(a-ik\om)\vphi_k+c\psi_k=0,\;\; \psi'_k+(d-ik\om)\psi_k+b\vphi_k=0
\ee
and the boundary conditions
\beq
\label{adjBC}
\psi_k(0)=r_0 \vphi_k(0),\; \vphi_k(1)=r_1\psi_k(1).
\ee
This yields, as in Section \ref{sec:iso}, that $(\vphi,\psi) \in \tilde{
V}^\ga(r_0,r_1)$ and
$(\tilde{A}+\tilde{B})(\vphi,\psi)=0$.
\end{proof}

Lemma \ref{lem:adj} implies that assertion (iii) of Theorem \ref{thm:fredh} 
is true. 
Hence, it remains to prove (\ref{5.2}), i.e.
\beq
\label{5.7}
\dim \ker (A+B)=\dim \ker(\tilde{A}+\tilde{B}).
\ee

\begin{lemma}\label{lem:finite}
There exists $k_0 \in \N$ such that for all $(u,v) \in \ker(A+B)$ and 
all $(\tilde{u},\tilde{v})
\in \ker(\tilde{A}+\tilde{B})$ and all $k \in \Z$ 
with $|k|>k_0$ we have $u_k=v_k=\tilde u_k=
\tilde v_k=0$.
\end{lemma}

\begin{proof} Suppose, contrary to our claim. Then there exists,
for example, a sequence
$(u^j,v^j)_{j \in \N}\in\ker(A+B)$ such that for all $j \in \N$ there 
is $k_j \in \Z$
with $|k_j| \ge j$ and $u^j_{k_j} \not=0$ or $v^j_{k_j}
\not=0$.
Without loss of generality we can assume that $k_j\not=k_l$ for $j\not=l$. 
Using the notation (\ref{e}) again, we see that the functions 
$(u^j_{k_j}e_{k_j},
v^j_{k_j}e_{k_j})$ belong to $\ker(A+B)$ and are linearly 
independent. On the other side we know that
$\ker(A+B)<\infty$, and this is a 
contradiction.
\end{proof}

Lemma \ref{lem:finite} implies that
$$
\ker(A+B)=\left\{\sum\limits_{|k|\le k_0} (u_ke_k,v_ke_k): \;(u_k,v_k) 
\mbox{ solves (\ref{eq:3.3}), (\ref{eq:3.4}) with } f_k=g_k=0 \right\},
$$
$$
\ker(\tilde{A}+\tilde{B})=\left\{\sum\limits_{|k|\le 
k_0} (\vphi_ke_k,\psi_ke_k): \;(\vphi_k,\psi_k) 
\mbox{ solves (\ref{adjGl}), (\ref{adjBC})}\right\}.
$$
It is known that, given $k \in \Z$, the number of linearly independent 
solutions to 
(\ref{eq:3.3}), (\ref{eq:3.4}) with  $f_k=g_k=0$ (this number is zero,
one or two) equals 
to the number of  linearly independent solutions to 
(\ref{adjGl}), (\ref{adjBC}). Hence, (\ref{5.7}) is proved.


\section{Closing remarks and open questions}\label{sec:clorem}
\setcounter{equation}{0}

In this final section we formulate some closing remarks, generalizations and open questions 
related to Theorem \ref{thm:fredh}.

\begin{rem}\label{rem:fredh} {\bf about \boldmath $L^\infty$ \unboldmath perturbations of  \boldmath $b$ \unboldmath and   
\boldmath $c$\unboldmath}
It seems to be an open question if the assumption $b,c \in BV(0,1)$ of 
Theorem \ref{thm:fredh}
can be weakened to $b,c \in L^\infty(0,1)$. But at least for ``almost 
all'' $b,c \in L^\infty(0,1)$
Theorem \ref{thm:fredh} remains to be true. More exactly, the following 
generalization of assertion (ii) of 
Theorem \ref{thm:fredh} holds:

Let $\ga \ge 1$, $a,d\in\L^{\infty}(0,1)$, and suppose
(\ref{eq:1.4}). Then there exists an open and dense set $M \subseteq (L^\infty(0,
1))^2$
such that $(BV(0,1))^2 \subset M$ and that for all $(b,c) \in M$
the operator $A+B$ is Fredholm of index zero
from $V^{\gamma}(r_0,r_1)$ into $W^{\gamma}$.

This generalization is true because the set of index zero Fredholm operators 
is open in
${\cal L}(V^{\gamma}(r_0,r_1);W^{\gamma})$, because the operator $B$ depends 
continuously (in the operator norm 
in ${\cal L}(W^{\gamma})$) on the coefficient functions $b$ and $c$ (in 
the $L^\infty$ norm), 
and because $BV(0,1)$ is dense in $L^\infty(0,1)$.
\end{rem}

\begin{rem}\label{rem:timedep}  {\bf about time depending perturbations of  \boldmath$a$\unboldmath,
\boldmath$b$\unboldmath,   \boldmath$c$ \unboldmath  and   \boldmath$d$\unboldmath}
The question, if an analog to Theorem \ref{thm:fredh} is true for general 
$2\pi/\om$-periodically 
time-depending coefficients $a$, $b$, $c$, and $d$ seems to be much more 
complicated (even if
$a$, $b$, $c$, and $d$ are supposed to be smooth). But at least for ``weakly 
time-depending'' coefficients 
Theorem \ref{thm:fredh} remains to be true. More exactly, the following holds:

Let $\ga \ge 1$, $a,d\in\L^{\infty}(0,1)$,  $b,c\in BV(0,1)$, and suppose
(\ref{eq:1.4}). Then there exists $\varepsilon>0$ such that the following 
is true: 
Take smooth functions $\tilde{a},\tilde{b},\tilde{
c},\tilde{d}: [0,1]\times\R  \to \R$ which are 
$2\pi/\om$-periodic with respect to the second argument. Suppose the $L^\infty$ 
norms of 
$\tilde{a},\tilde{b},\tilde{c}$, and $\tilde{
d}$ to be less than $\varepsilon$. Define the operators
$A \in {\cal L}(V^{\gamma}(r_0,r_1);W^{\gamma})$ and $B \in {\cal L}(W^{\gamma})$ 
as above by 
replacing $a$, $b$, $c$, and $d$ by $a+\tilde{a}$, $b+\tilde{
b}$, $c+\tilde{c}$, and $d+\tilde{d}$, 
respectively. Then $A$ is an isomorphism from $V^{\gamma}(r_0,r_1)$ onto 
$W^{\gamma}$, and $A+B$ is
Fredholm of index zero from $V^{\gamma}(r_0,r_1)$ into~$W^{\gamma}$.

The argument of the proof of this assertion is, again, the openness of the sets 
of isomorphisms and of index zero Fredholm operators and the continuous 
dependence of the 
operators $A$ and $B$ on the functions $\tilde{a},\tilde{
b},\tilde{c}$, and $\tilde{d}$.
\end{rem}

\subsection*{Acknowledgments}
This work was done while 
the first author
visited the Institute of Mathematics of the
Humboldt University of Berlin.
She is thankful to Prof. Hans-J{\"u}rgen~Pr{\"o}mel for his kind
hospitality during her stay at the Humboldt University.

The authors would like to thank the anonymous referee for his  
comments that lead to improvements of the paper.

\end{document}